\documentclass{aspm}

\articleinfo{XXX}{YYY}{ZZZ}
\usepackage{amsmath,amssymb}

\usepackage[letterpaper]{geometry}

\newcommand{\Span}[1]{\left<#1\right>}
\newcommand{\broken}{\dashrightarrow}

\newcommand{\onto}{\twoheadrightarrow}
\newcommand{\iso}{\simeq}

\newcommand{\sA}{\mathcal A}
\newcommand{\sB}{\mathcal B}
\newcommand{\sI}{\mathcal I}
\newcommand{\sL}{\mathcal L}
\newcommand{\sM}{\mathcal M}
\newcommand{\Oh}{\mathcal O}

\newcommand{\C}{\mathbb C}
\newcommand{\FF}{\mathbb F}
\newcommand{\PP}{\mathbb P}
\newcommand{\Q}{\mathbb Q}
\newcommand{\Z}{\mathbb Z}

\newcommand{\Chat}{\widehat C}
\newcommand{\Dhat}{\widehat D}
\newcommand{\Ehat}{\widehat E}
\newcommand{\Lhat}{\widehat L}
\newcommand{\Mhat}{\widehat M}
\newcommand{\Phat}{\widehat P}
\newcommand{\Xhat}{\widehat X}
\newcommand{\Lahat}{\widehat\La}
\newcommand{\qhat}{\widehat q}
\newcommand{\sBhat}{\widehat{\sB}}
\newcommand{\sLhat}{\widehat{\sL}}
\newcommand{\sMhat}{\widehat{\sM}}
\newcommand{\Xihat}{\widehat{\Xi}}
\newcommand{\Thhat}{\widehat{\Theta}}

\newcommand{\Abar}{\overline A}
\newcommand{\Cbar}{\overline C}
\newcommand{\Dbar}{\overline D}
\newcommand{\Ebar}{\overline E}
\newcommand{\Fbar}{\overline F}
\newcommand{\Gbar}{\overline G}
\newcommand{\Lbar}{\overline L}
\newcommand{\Pbar}{\overline P}

\newcommand{\Xbar}{\overline X}
\newcommand{\Ybar}{\overline Y}
\newcommand{\fbar}{\overline f}
\newcommand{\Labar}{\overline\Lambda}
\newcommand{\sLbar}{\overline{\sL}}
\newcommand{\sMbar}{\overline{\sM}}
\newcommand{\Xibar}{\overline{\Xi}}
\newcommand{\wD}{\widetilde D}
\newcommand{\wE}{\widetilde E}
\newcommand{\wF}{\widetilde F}

\newcommand{\wX}{\widetilde X}

\newcommand{\wXi}{\widetilde{\Xi}}
\newcommand{\wM}{\widetilde{\sM}}
\newcommand{\al}{\alpha}
\newcommand{\be}{\beta}
\newcommand{\ga}{\gamma}
\newcommand{\Ga}{\Gamma}
\newcommand{\de}{\delta}
\newcommand{\De}{\Delta}
\newcommand{\la}{\lambda}
\newcommand{\La}{\Lambda}

\newcommand{\ze}{\zeta}

\DeclareMathOperator{\ct}{ct}
\DeclareMathOperator{\g}{g}
\DeclareMathOperator{\lcm}{lcm}
\DeclareMathOperator{\rk}{rk}
\DeclareMathOperator{\Bs}{Bs}
\DeclareMathOperator{\Cl}{Cl}
\DeclareMathOperator{\Gr}{Gr}

\DeclareMathOperator{\Pic}{Pic}

\DeclareMathOperator{\qW}{q_F}
\DeclareMathOperator{\qQ}{q_\Q}

\DeclareMathOperator{\Supp}{Supp}
\DeclareMathOperator{\Tors}{Tors}
\DeclareMathOperator{\mult}{mult}

\newcommand{\muu}{{\boldsymbol{\mu}}}
\newcommand{\qq}{\mathbin{\sim_{\scriptscriptstyle{\Q}}}}

\numberwithin{equation}{subsection}

\newtheorem{lem}[subsection]{Lemma}
\newtheorem{prop}[subsection]{Proposition}

\newtheorem{stheorem}[equation]{}

\theoremstyle{definition}
\newtheorem{rem}[equation]{Remark}
\newtheorem{exa}[equation]{Example}

\newtheorem{case}[subsection]{}
\newtheorem{scase}[equation]{}

\title{On $\Q$-Fano $3$-folds of Fano index 2}
\dedicatory{To Shigefumi Mori in friendship and admiration}
\author{Yuri Prokhorov \and Miles Reid}
\thanks{Y.P. acknowledges partial support from RFBR grants No.\
11-01-92613-KO\_a, 11-01-00336-a, the grant of Leading Scientific
Schools No.\ 2998.2014.1, and Russian Academic Excellence Project
`5--100'. \endgraf
M.R. is partly supported as Scholar of KIAS}
\address{Yuri Prokhorov, Steklov Mathematical Institute,
\newline 8 Gubkina street, Moscow 119991, Russia
\newline\indent
Department of Algebra, Faculty of Mathematics,
\newline Moscow State University, Moscow 117234, Russia
\newline\indent
Laboratory of Algebraic Geometry, \newline
National Research University Higher School of Economics, Russia}
\email{prokhoro@mi.ras.ru}

\address{Miles Reid, Mathematics Institute, University of Warwick,
\newline Coventry CV4 7AL, UK}
\email{Miles.Reid@warwick.ac.uk}

\begin{document}
\setcounter{page}{397}
\maketitle


\section{Introduction}
\subsection{$\Q$-Fano $3$-folds}
A {\em $\Q$-Fano $3$-fold} is a projective $3$-fold $X$ with at worst
terminal singularities and ample anticanonical divisor $-K_X$. Here,
bearing in mind Mori's fundamental notion of extremal ray, we assume
also that $X$ is $\Q$-factorial and has rank~1, that is, $\Pic X\iso \Z$
or equivalently, $\Cl X\otimes\Q\iso \Q$. We define the {\em Fano} and
{\em $\Q$-Fano index} of $X$ by:
\begin{align*}
\qW(X)&:= \max \{ q\in\Z \mid \hbox{$-K_X\sim qA$ with $A$ a Weil divisor}\},
\\[3pt]
\qQ(X)&:= \max \{ q\in \Z \mid \hbox{$-K_X\qq qA$ with $A$ a Weil divisor}\},
\end{align*}
where $\sim$ is linear equivalence and $\qq$ is $\Q$-linear equivalence.
Clearly, $\qW(X)$ divides $\qQ(X)$, and the two coincide unless
$K_X+qA\in\Cl X$ is a nontrivial torsion element. An important invariant
of a $\Q$-Fano $3$-fold is its {\em genus} $\g(X):=\dim|{-}K_X|-1$.

\subsection{Background facts}
Kaori Suzuki \cite{Suzuki-2004} restricts the $\Q$-Fano index of $X$ to
one of
\begin{equation}
\label{eq!Kaori}
\qQ(X)\in \{1,\dots, 11, 13,17,19\}.
\end{equation}
See also \cite[Lemma~3.3]{Pr2008a}. Moreover, the following
results are due to the first author.

\begin{stheorem}{\bf Theorem (\cite{Pr2008a}).}
\label{th!bk1}
Let $X$ be a $\Q$-Fano $3$-fold of\/ $\Q$-Fano index $q:=\qQ(X)\ge9$.
Then $\Cl X \iso \Z$.
\begin{enumerate}
\renewcommand{\labelenumi}{(\roman{enumi})}
\renewcommand\theenumi{(\roman{enumi})}

\item\label{th!bk-q=19}
If $q=19$ then $X\iso \PP(3,4,5,7)$.

\item\label{th!bk-q=17}
If $q=17$ then $X\iso \PP(2,3,5,7)$.

\item\label{th!bk-q=13}
If $q=13$ and $\g(X)>4$ then $X\iso \PP(1,3,4,5)$.

\item\label{th!bk-q=11}
If $q=11$ and $\g(X) >10$ then $X\iso \PP(1,2,3,5)$.

\item\label{th!bk-q=10}
$q\ne10$.
\end{enumerate}
\end{stheorem}

\begin{stheorem}{\bf Theorem (\cite{Pr2010-i}).}
\label{th!bk2}
Let $X$ be a $\Q$-Fano $3$-fold of\/ $\Q$-Fano index $q$.

\begin{enumerate}
\renewcommand{\labelenumi}{(\roman{enumi})}
\renewcommand\theenumi{(\roman{enumi})}
\setcounter{enumi}{5}
\item\label{th!bk-q=9}
If $q=9$ and $\g(X)>4$ then
$X\iso X_6 \subset \PP(1,2,3,4,5)$.

\item\label{th!bk-q=8}
If $q=8$ and $\g(X)>10$ then
$X\iso X_6\subset \PP(1,2,3^2,5)$ or
$X_{10}\subset\PP(1,2,3,5,7)$.

\item\label{th!bk-q=7}
If $q=7$ and $\g(X)>17$ then
$X\iso \PP(1^2,2,3)$.

\item\label{th!bk-q=6}
If $q=6$ and $\g(X)>15$ then
$X\iso X_6\subset \PP(1^2,2,3,5)$.

\item\label{th!bk-q=5}
If $q=5$ and $\g(X)>18$ then $X\iso \PP(1^3,2)$ or
$X_4\subset \PP(1^2,2^2,3)$.

\item\label{th!bk-q=4}
If $q=4$ and $\g(X)>21$ then
$X\iso \PP^3$ or $X_4\subset \PP(1^3,2,3)$.

\item\label{th!bk-q=3}
If $q=3$ and $\g(X)>20$ then
$X\iso X_2\subset \PP^4$ or $X_3\subset \PP(1^4,2)$.
\end{enumerate}
\end{stheorem}
Here we study the case $\qQ(X)=2$.

\begin{stheorem}{\bf Theorem (\cite{BS07mm}).}
The Hilbert series of\/ $\Q$-Fano $3$-folds with $q=\qQ(X)=\qW(X)=2$
belong to at most $1492$ cases.

\end{stheorem}

The online database \cite{GRD} lists the numerical type of candidates
(the data going into the Hilbert series of their graded rings).

\subsection{Main results}

\begin{stheorem}{\bf Main Theorem.}\label{th!m}
Let $X$ be a $\Q$-Fano $3$-fold of rank $1$ with $\qQ(X)=\qW(X)=2$ and
$K_X$ not Cartier. Let $A$ be a Weil divisor on $X$ such that $-K_X=2A$.

Then $\dim|A|\le 4$. Moreover, if\/ $\dim|A|=4$, then $X$ belongs to the
single irreducible family constructed in \ref{ex!main} (see also
\ref{ss!MEx}).
\end{stheorem}

\begin{stheorem}{\bf Corollary.}
A $\Q$-Fano $3$-fold with $\qQ(X)=\qW(X)=2$ and $K_X$ not Cartier has
$\g(X)\le 16$.
\end{stheorem}

\begin{rem} If $K_X$ is Cartier and $\qW(X)=2$, then $X$ is a {\em del
Pezzo variety} \cite{Fujita-1162108}. Two cases with $\Cl X\iso \Z$
have $\dim|A|>4$:
\begin{enumerate}
\renewcommand{\labelenumi}{(\alph{enumi})}
\renewcommand\theenumi{(\alph{enumi})}

\item the complete intersection of two quadrics $X=X_{2\cdot 2}\subset
\PP^5$, with $\dim|A|=5$ and $\g(X)=19$; and

\item $X=X_5\subset\PP^6$ a section of the Grassmannian
$\Gr(2,5)\subset\PP^9$ by a subspace of codimension $3$, with
$\dim|A|=6$ and $\g(X)=23$. In this case $X$ must be smooth by
\cite[Cor.~5.3]{Pr2010a}.
\end{enumerate}
\end{rem}

\subsection{Strategy of proof}
Sections~\ref{s!pf1}--\ref{s!pf2} contain the proof of Main
Theorem~\ref{th!m}. The Kawamata blowup of a $\frac1r(2,a,r-a)$ point
initiates a Sarkisov link ending in a fibre space or a $\Q$-Fano
$3$-fold with $q\ge3$; the assumption $\dim|A|\ge4$ leads to a
manageable case division. The auxiliary Section~\ref{s!pf0} treats the
cases with $q\ge3$, most of which lead to a contradiction, with just one
surviving in Section~\ref{s!pf2} to characterize our Main Example.

\subsection{The Main Example} \label{ss!MEx}
Section~\ref{s!mEx} gives several constructions of the exceptional
family of Main Theorem~\ref{th!m}, $\Q$-Fanos $X$ with $\Cl X=\Z\cdot
A$, $K_X=-2A$ and $\dim|A|=4$. They arise from the simplest type of
Sarkisov link:
\begin{equation}
\begin{matrix}
\kern-1cm E \,\subset\, X_1 \,\supset\, F \kern-1cm \\[2pt]
\kern-1cm \begin{picture}(0,0)(0,0)
\put(-7,9){\vector(-1,-2){7}}
\put(6,9){\vector(1,-2){7}}
\end{picture}
\kern-1cm \\[3pt]
\kern-5.3em \hphantom{\owns P_0} P \,\in\, X \kern2.6em Q \,\supset\, \Ga_5 \,\owns\, P_0 \kern-6em
\end{matrix}
\label{eq!Sark}
\end{equation}
starting from the nonsingular quadric hypersurface $Q\subset\PP^4$, a
point $P_0\in Q$ and an irreducible curve $\Ga_5\subset E_0$ of degree~5
contained in the tangent hyperplane section $E_0=Q\cap T_{P_0}Q$, with
$\mult_{P_0}\Ga_5=3$. We make the symbolic blowup $X_1\to Q$ of $\Ga_5$,
then contract the birational transform $E_1\iso E_0\iso\PP(2,1,1)\subset
X_1$ of $E_0$ to a $\frac13(1,2,2)$ orbifold point.

The symbolic blowup of $\Ga\subset E_0\subset Q$ is the relative Proj of
the symbolic algebra $\sA=\bigoplus \sI_\Ga^{[n]}$. For a singular curve
$\Ga$ contained in a nodal surface, this is a local graded ring
construction with a universal description, studied in much more detail
and generality in Tom Ducat's thesis \cite{Du15}. Compare \cite{Du14}.

\subsection{Discussion}
The study of $\Q$-Fanos divides into birational and biregular
considerations. Biregular methods study projective embedding by
multiples of $A$, or more precisely, generators and relations for the
Gorenstein graded ring $R(X,A)$. This is effective when $R(X,A)$ has
small codimension, especially if it is a hypersurface or codimension~2
complete intersection, etc. In contrast, birational methods are powerful
when the linear system $|A|$ is large, implying a low canonical
threshold, and allowing us to impose noncanonical singularities on $|A|$
and study $X$ via the resulting Sarkisov link, aiming for a birational
construction or a nonexistence result. The interest of this paper is as
a meeting point of the two methods.

\subsection{The fabulous half-elephant; more cases with $q=2$}
A surface section $F\in|A|$ of a $\Q$-Fano $3$-fold $X$ of index~2 is a
del Pezzo surface (sometimes very singular). In a few cases where $F$ has
the simplest orbifold points such as $\frac13(2,2)$ or $\frac15(2,4)$,
Reid and Suzuki \cite{RS03} study such surfaces in terms of cascades of
projections from nonsingular points. This foreshadows one construction
of our Main Example in Section~\ref{s!mEx}, and hints at other cases
that might make interesting challenges, especially the $X$ with
$\dim|A|=3$ or $2$. Del Pezzo surfaces with only $\frac13(2,2)$ orbifold
points are classified in current work of Alessio Corti and Liana
Heuberger \cite{CH15}. Kuzma Khrabrov \cite{Kh14} has partial results on
$\Q$-Fano $3$-folds $X$ of index~2 with $\dim|A|\ge2$.

\section{The method} \label{sect-2}

\subsection{Construction of a Sarkisov link
\cite{Alexeev-1994ge}}\label{ss!Slink}
Let $\sM$ be a linear system on $X$ with no fixed part, and canonical
threshold $c:=\ct(X,\sM)$. Assume $-(K_X+c\sM)$ is ample. Then
$(X,c\sM)$ is canonical but not terminal, so we can pull out an
irreducible divisor $E\subset\wX$ by an extremal divisorial extraction
$f\colon \wX\to X$, such that $\wX$ has only terminal $\Q$-factorial
singularities, $\rho(\wX/X)=1$, and $f$ is $(K+c\sM)$-crepant:
\begin{equation}
\label{eq!KX+cM}
K_{\wX}+c\wM =f^*(K_X+c\sM).
\end{equation}
As in \cite{Alexeev-1994ge}, running a $(K+c\sM)$-MMP on $\wX$ gives a
Sarkisov link of type~I or~II:
\begin{equation}
\renewcommand{\arraycolsep}{.2em}
\label{eq3.1}
\begin{matrix}
&& \wX & \broken & \Xbar \\[2pt]
& \begin{picture}(0,0)(0,0)
\put(-9,3){$f$}
\put(3,9){\vector(-1,-2){7}}
\end{picture} &&&
\begin{picture}(0,0)(0,0)
\put(12,3){$\fbar$}
\put(5,9){\vector(1,-2){7}}
\end{picture} \\[3pt]
X &&&&&& \Xhat
\end{matrix}
\end{equation}
where $\wX$ and $\Xbar$ have only $\Q$-factorial terminal singularities,
$\rho(\wX)=\rho(\Xbar)=2$, $\wX\broken\Xbar$ is a chain of log flips,
and $\fbar$ is a Mori extremal contraction, either a divisorial
contraction to a $\Q$-Fano $3$-fold $\Xhat$, or a Mori fibre space over
a curve or surface $\Xhat$. In either case, $\rho(\Xhat)=1$. We write
$\wD$ and $\Dbar$ for the birational transform on $\wX$ and $\Xbar$ of a
divisor or linear system $D$ on $X$.

Assume that $K_X+\la\sM+\Xi\qq0$ for some $\la>c$ and an effective
$\Q$-divisor $\Xi$. We can write
\begin{equation}
\label{eq!KX+laM}
K_{\wX}+\la \wM + \wXi+aE\qq f^*(K_X+\la \sM+\Xi)\qq0,
\end{equation}
where $a>0$ is the log discrepancy of $f$. Note that if
$K_X+\la\sM+\Xi\sim 0$ then it is a Cartier divisor; we can assume that
the mobile system $\wM$ has no common divisor with $\wXi$. Then
$K_{\wX}$, $\la\wM$ and $\wXi$ are all integral Weil divisors, and
therefore $a$ is an integer.

\begin{rem}\label{rem!E}
We use the extremal extraction $\wX\to X$ with ray $R$ and exceptional
surface $E$ to initiate a Sarkisov link. By \eqref{eq!KX+cM},
$K_{\wX}+c\wM$ is nef and big, with
\begin{equation}
K_{\wX}\cdot R < 0, \quad (K_{\wX}+c\wM)\cdot R = 0, \quad\hbox{so that
$\wM\cdot R>0$}.
\end{equation}
The MMP that constructs the Sarkisov link proceeds by increasing $\la$
in $K+\la\sM$. Each step makes $K+\la\sM$ bigger on the ray $R$, so on
the exceptional surface $E$ and its birational transforms. Thus the MMP
can never contract the birational transform of $E$.
\end{rem}

\subsection{Case $\fbar$ not birational}
\label{ss!Mfs}
Assume that $\fbar$ is not birational. Then $\Xhat$ is either a smooth
rational curve or a del Pezzo surface with at worst Du Val singularities
and $\rho(\Xhat)=1$ \cite{Mori-Prokhorov-2008}. We also have
$\fbar(\Ebar)=\Xhat$ by Remark~\ref{rem!E}, or because no multiple
$n\Ebar$ of the exceptional divisor $\Ebar$ of $f$ moves on $\Xbar$. In
this case we write $\Fbar$ for a general fiber of $\fbar$. Let $\Theta$
be an ample Weil divisor on $\Xhat$ whose class generates $\Cl \Xhat$
modulo torsion. If $\Xhat$ is a surface with $K_{\Xhat}^2=1$, we take
$\Theta=-K_{\Xhat}$.

\begin{scase}\label{possibilities-Du-Val-del-Pezzo}
For $\Xhat$ a surface, one of the following holds:
\begin{enumerate}
\renewcommand{\labelenumi}{(\roman{enumi})}
\renewcommand\theenumi{(\roman{enumi})}

\item $-K_{\Xhat}\cdot \Theta=3$, $-K_{\Xhat}\sim 3\Theta$, $\Xhat\iso
\PP^2$ and $\dim|\Theta|=2$;

\item $-K_{\Xhat}\cdot \Theta=2$, $-K_{\Xhat}\sim 4\Theta$, $\Xhat\iso
\PP(1,1,2)$ and $\dim|\Theta|=1$;

\item $-K_{\Xhat}\cdot \Theta=1$, $-K_{\Xhat}\sim d\Theta$, where
$d:=K_{\Xhat}^2\le 6$, and the minimal resolution of $\Xhat$ is a blowup
of $\PP^2$ at $9-d$ points in {\em almost general position}. In this
case, $\dim|\Theta|=0$ or $1$. Moreover, by Kawamata--Viehweg vanishing
and orbifold Riemann--Roch \cite{YPG}, for an ample Weil divisor $B\qq
t\Theta$ we have
\begin{equation}
\dim|B|\le\frac{t(t+d)}{2d}.
\end{equation}
\end{enumerate}
\end{scase}

\subsection{Case $\fbar$ birational}
Assume that the contraction $\fbar$ is birational. In this case, $\Xhat$
is a $\Q$-Fano $3$-fold and $\fbar$ contracts a unique exceptional
divisor $\Fbar$. Remark~\ref{rem!E} implies that $\Ebar\ne\Fbar$ (or
argue that $\Ebar=\Fbar$ would imply that $X\broken\Xhat$ is an
isomorphism in codimension, leading to a contradiction). Write
$\wF\subset\wX$ and $F:=f(\wF)$ for its birational transform. Set
$\qhat:=\qQ(\Xhat)$. For a divisor $\Dbar$ on $\Xbar$, we put
$\Dhat:=\fbar_*\Dbar$.

\subsection{Computer search for $\Q$-Fano $3$-folds}
\label{ss!al}
All $\Q$-Fano $3$-folds belong to a finite number of algebraic families
\cite{Ka92bF}. In fact, Kawamata's proof implies that the possible
``candidate'' $\Q$-Fano $3$-folds can be listed, although the volume of
computation makes computer searches inevitable. This method was used in
\cite{Suzuki-2004}, \cite{BS07j}, \cite{BS07mm},
\cite{Pr2007-Qe}, \cite{Pr2008a},
\cite{Pr2010-i}. See \cite{GRD} for explicit lists.

We outline the algorithm, starting with a useful remark.

\begin{rem} \label{remark-index-odd}
The local analytic Weil divisor class group of a \hbox{$3$-fold}
$\Q$-factorial terminal point $P\in X$ is cyclic $\Cl(X,P)\iso \Z/r$,
and is generated by the canonical divisor $K_X$
\cite[Lemma~5.1]{Kawamata-1988-crep}. In particular, if $X$ is a
$\Q$-Fano $3$-fold, its local Gorenstein index $r$ at every terminal
point is coprime to the $\Q$-Fano index $q=\qW(X)$.
\end{rem}

\begin{scase}
Let $X$ be a $\Q$-Fano $3$-fold. For simplicity we assume that
$q:=\qQ(X)=\qW(X)\ge3$ (the only case we need in this section). Let $A$
be a Weil divisor such that $-K_X\sim qA$ and $\sB(X)=\{(r_P,b_P)\}$
the basket of orbifold points of $X$ \cite{YPG}.
\end{scase}

\textbf{Step 1.}
We have the equality
\begin{equation}
-K_X\cdot c_2(X)+\sum_{P\in \sB} \frac{r_P-1}{r_P}= 24,
\end{equation}
where $-K_X\cdot c_2(X)>0$ \cite{Ka92bF}. Hence there is only a
finite (but huge) number of possibilities for the basket $\sB(X)$ and
$-K_X\cdot c_2(X)$. Let $r:=\lcm(\{r_P\})$ be the Gorenstein index of
$X$.

\medskip
\textbf{Step 2.}
\eqref{eq!Kaori} says that $q\in \{3,\dots,11, 13,17,19\}$.
Remark~\ref{remark-index-odd} implies that $\gcd(q,r)=1$, which
eliminates some possibilities.

\medskip
\textbf{Step 3.}
In each case we compute $A^3$ by the formula
\begin{equation*}
A^3=\frac{12}{(q-1)(q-2)}\Bigl(
1-\frac{A\cdot c_2}{12}+\sum_{P\in B} c_P(-A)
\Bigr).
\end{equation*}
(see \cite{Suzuki-2004}), where $c_P$ is the correction term in the
orbifold Riemann--Roch formula \cite{YPG}. The number $rA^3$
must be an integer \cite[Lemma~1.2]{Suzuki-2004}.

\medskip
\textbf{Step 4.}
Next, the Bogomolov--Miyaoka inequality (see \cite{Ka92bF})
implies that
\begin{equation}
\left(4q^2 - 3q\right)A^3\le -4K_X\cdot c_2(X)
\end{equation}
\cite[Prop.\ 2.2]{Suzuki-2004}.

\medskip
\textbf{Step 5.}
Finally, the Kawamata--Viehweg vanishing theorem gives
$\chi(tA)=h^0(tA)=0$ for $-q<t<0$. We check this condition using
orbifold Riemann--Roch \cite{YPG}, \cite{Ice}.

\section{On $\Q$-Fano $3$-folds of Fano index $\ge3$}\label{s!pf0}
\subsection{A result of Fujita}
A {\em polarized variety} is a pair $(X,S)$ consisting of a projective
variety $X$ and an ample Cartier divisor $S$ on $X$. Its {\em
$\De$-genus} is defined as follows \cite{Fujita-1162108}:
\begin{equation}
\De(X,S)=\dim X+S^{\dim X}-h^0(X,\Oh_X(S)).
\end{equation}
It is known that $\De(X,S)\ge0$ and Fujita \cite{Fujita-1162108}
classifies polarized varieties of small $\De$-genera. We use the
following easy consequence of Fujita's classification.

\begin{stheorem}{\bf Lemma.}\label{lemma-del-pezzo}
Let $X$ be a $\Q$-Fano $3$-fold and $S$ an ample Weil divisor on $X$
such that $\dim|S|>0$, $|S|$ has no fixed components, and $-K_X\qq \la
S$ with $\la \ge2$. Assume that the pair $(X,|S|)$ is terminal. Then one
of the following holds:
\begin{enumerate}
\renewcommand{\labelenumi}{(\roman{enumi})}
\renewcommand\theenumi{(\roman{enumi})}

\item $X\iso \PP^3$, $\la =4$, $\dim|S|=3$;

\item $X\iso \PP^3$, $\la =2$, $\dim|S|=9$;

\item $X\iso X_2\subset \PP^4$ is a smooth quadric, $\la =3$,
$\dim|S|=4$;

\item $X$ is a del Pezzo $3$-fold of degree $1\le d\le 5$, $\la =2$,
$\dim|S|=d+1$;

\item $X\iso \PP(1^3,2)$, $\la =5/2$, $\dim|S|=6$.
\end{enumerate}
\end{stheorem}

\begin{proof}
Replace $S$ with a general member of $|S|$. Since $(X,|S|)$ is terminal,
the surface $S$ is smooth and contained in the smooth locus of $X$
\cite[1.22]{Alexeev-1994ge}. By the adjunction formula we have $-K_S\sim
(\la -1)S_{|S}$. Hence $S$ is a (smooth) del Pezzo surface and $(\la -1)^2
S^3=K_S^2$. Since $H^i(X,\Oh_X)=0$ and $H^i(S,\Oh_S(S))=0$ for $i>0$, by
Riemann--Roch we have
\begin{equation}
h^0(X,\Oh_X(S))=h^0(S,\Oh_S(S))+1=
\frac \la 2 S^3+2.
\end{equation}
Therefore,
\begin{equation}
\De(X,S)=3+ S^3-\frac \la 2 S^3-2=1+\frac{(2-\la )S^3}2=
1+\frac{(2-\la )K_S^2}{2(\la -1)^2}.
\end{equation}
If $S\iso\PP^2$, then $\Oh_S(S)=\Oh_{\PP^2}(l)$, where $3=(\la
-1)l\ge l$. Then $\De(X,S)=0$ and \cite[Th.\ 5.10 and
5.15]{Fujita-1162108} gives cases (i) and (v). If $S\iso\PP^1\times
\PP^1$, then $\Oh_S(S)=\Oh_{\PP^1\times \PP^1}(k,k)$, where
$k(\la-1)=2$. So, $\la=2$ or $3$, $\De(X,S)=0$, and
\cite[Th.\ 5.10 and 5.15]{Fujita-1162108} gives cases (ii) or (iii).
Finally, if $S\not \iso\PP^2$, $\PP^1\times \PP^1$, then $K_S$ is a
primitive element of $\Pic S$. Hence $\la=2$ and $\De(X,S)=1$.
Then we have case (iv) \cite[Ch.\ 1, \S 9]{Fujita-1162108}.
\end{proof}

\begin{lem}[{\cite[Th.\ 1.4 (vii)]{Pr2008a}}] \label{Lemma-q=5}
Let $X$ be a $\Q$-Fano $3$-fold with terminal singularities and with
$q:=\qQ(X)\ge5$. Let $A$ be a Weil divisor such that $-K_X\qq qA$.
If\/ $\dim|A|\ge2$, then $X\iso \PP(1^3,2)$.
\end{lem}

\begin{proof}
We first consider the case $\rk\Cl X=1$ and $\qQ(X)=\qW(X)$ (in
particular, $X$ is a $\Q$-Fano $3$-fold and $-K_X\sim q A$). Running a
computer search as in \ref{ss!al} gives $-K_X^3\ge125/2$. Then by
\cite{Pr2007-Qe} we have $X\iso \PP(1^3,2)$.

Next consider the case $\rk\Cl X>1$ and $\qQ(X)=\qW(X)$. We get a
contradiction in this case. Run the MMP. The property $-K_X\sim q A$ is
preserved. At the end we get a $\Q$-Fano $3$-fold $\Xbar$ with
$K_{\Xbar}\sim q \Abar$ and $\dim|\Abar|\ge2$. By \eqref{eq!Kaori} we
have $\qQ(\Xbar)=\qW(\Xbar)=q$. By the above $\Xbar\iso \PP(1^3,2)$
and $\dim|A|=\dim|\Abar|=2$. Let $\Pbar\in \Xbar$ be the point of type
$\frac12(1,1,1)$. Consider the final step $g\colon \wX\to \Xbar$ of the
MMP, a divisorial contraction, and let $\wE\subset\wX$ be its
exceptional divisor. There are the following possibilities:
\begin{enumerate}
\renewcommand{\labelenumi}{(\alph{enumi})}
\renewcommand\theenumi{(\alph{enumi})}

\item
$g(\wE)=\Pbar$. Then $K_{\wX}\qq g^*K_{\Xbar}+\frac12 \wE$,
$\wE\iso \PP^2$, and $\Oh_{\wE}(\wE)\iso
\Oh_{\PP^2}(-2)$ \cite{Kawamata-1996}. Hence, $\Oh_{\tilde
E}(-K_{\wX})\iso \Oh_{\PP^2}(1)$. We get a contradiction because
$-K_{\wX}$ is divisible by $q\ge5$.

\item $g(\wE)$ is either a smooth point or a curve. In this case
$g(\wE)\not\subset \Bs |\Abar|=\{\Pbar\}$. On the other hand, $g$ is a
$K_{\wX}$-negative contraction, a contradiction.
\end{enumerate}

Finally assume that the torsion part of $\Cl X$ is nontrivial. Every
torsion element $\xi_1\in \Cl X$ of order $n_1>1$ defines a
$\muu_{n_1}$-cover $\pi_1\colon X_1\to X$ that is \'etale in
codimension~2. Repeating the procedure we get a sequence
\begin{equation}
X_m \xrightarrow{\pi_m}
 X_{m-1} \xrightarrow{\pi_{m-1}}
\cdots \xrightarrow{\pi_2} X_1 \xrightarrow{\pi_1} X.
\end{equation}
with each $\pi_k$ a $\muu_{n_k}$-cover that is \'etale in codimension~2
and $\Cl X_m$ torsion free. By the above, $X_m\iso \PP(1^3,2)$. Since
\begin{equation}
h^0(X_m,\pi^*A)=h^0(X,A)=3,
\end{equation}
$\muu_{n_m}$ acts trivially on
$H^0(X_m,\pi^*A)=H^0(\Oh_{\PP(1^3,2)}(1))$. On the other hand, we can
take independent sections $x_1,x_2, x_3\in H^0(\Oh_{\PP(1^3,2)}(1))$ as
orbinates at the $\frac12(1,1,1)$-point $P_m\in X_m$. This contradicts
that the point $(X_m,P_m)/\muu_{n_m}$ is terminal.
\end{proof}

In a similar way to Lemma~\ref{Lemma-q=5}, one can prove the following.

\begin{lem}[{\cite[Th.\ 1.4 (vi)]{Pr2008a}}]
\label{Lemma-q=7}
Let $X$ be a $\Q$-Fano $3$-fold with terminal singularities and with
$q:=\qQ(X)\ge7$. Let $A$ be a Weil divisor such that $-K_X\qq qA$. If
$\dim|A|\ge1$, then $X\iso \PP(1^2,2,3)$.
\end{lem}

\begin{prop} \label{p!dim-2M}
Let $X$ be a $\Q$-Fano $3$-fold and let $q:=\qQ(X)$. Let $\sM$ be a
linear system on $X$ such that $\dim \sM\ge4$ and $-K_X\sim 2 \sM+\Xi$,
where $\Xi$ is a nonzero effective Weil divisor. Then $\Cl X\iso
\Z\cdot\Xi$, $q=2n+1$ is odd, and $\sM\sim n\Xi$. Moreover, one of the
following holds:
\begin{enumerate}
\renewcommand{\labelenumi}{(\roman{enumi})}
\renewcommand\theenumi{(\roman{enumi})}

\item \label{p!dim-2M-q=13}
$q=13$,\quad $X\iso \PP(1,3, 4, 5)$;

\item \label{p!dim-2M-q=11}
$q=11$,\quad $X\iso \PP(1,2, 3, 5)$;

\item \label{p!dim-2M-q=9}
$q=9$,\quad $X\iso X_6\subset \PP(1, 2, 3,4,5)$;

\item \label{p!dim-2M-q=7}
$q=7$,\quad $X\iso \PP(1^2, 2, 3)$;

\item \label{p!dim-2M-q=5a}
$q=5$,\quad $X\iso X_4 \subset \PP(1^2,2^2,3)$;

\item \label{p!dim-2M-q=5}
$q=5$,\quad $X\iso \PP(1^3,2)$;

\item \label{p!dim-2M-q=3}
$q=3$,\quad $X\iso X_2 \subset \PP^4$.
\end{enumerate}
\end{prop}

\begin{proof}
By assumption $q\ge3$. If $q\ge9$, then the assertion follows by
\cite[Prop.\ 3.6]{Pr2008a} and Theorem~\ref{th!bk2} \ref{th!bk-q=9}. So
assume that $3\le q\le 8$.

Let $A$ be a Weil divisor with $-K_X\qq\,qA$ and $n$ the integer such
that $\sM\qq nA$. If $\Cl X$ is torsion free, we can run the computer
search \ref{ss!al}. We get $q\ne4$ and $\g(\Xhat)\ge21$. Then by
Theorem~\ref{th!bk2} we get one of cases
\ref{p!dim-2M-q=7}--\ref{p!dim-2M-q=3}. Thus from now on we assume that
$\Cl X$ contains a nontrivial torsion element.

We may assume that $\sM$ has no fixed part. If the pair $(X,\sM)$ is
terminal, then $X$ is in \ref{p!dim-2M-q=5} or \ref{p!dim-2M-q=3} by
Lemma~\ref{lemma-del-pezzo}. Assume that $(X,\sM)$ is not terminal.
Apply Construction~\ref{ss!Slink} to $(X,\sM)$. We can write
\begin{equation*}
K_{\wX}+2\wM +\wXi+a\wE\sim f^*(K_X+2\sM+\Xi)\sim 0,
\end{equation*}
where $a\in \Z_{>0}$. Hence,
\begin{equation}\label{eq!KXbar}
K_{\Xbar}+2\sMbar+\Xibar+a\Ebar\sim 0.
\end{equation}

First consider the case of \ref{ss!Mfs} where $\fbar$ is not birational.
In particular, $\Xhat$ is either $\PP^1$ or a del Pezzo surface as in
\ref{possibilities-Du-Val-del-Pezzo}.

Assume that $\sMbar$ is $\fbar$-horizontal. Restricting the relation
\eqref{eq!KXbar} to a general fiber $\Fbar$ of $\fbar$
we get
\begin{equation}
-K_{\Fbar}\qq 2 \sMbar_{|\Fbar}+\Xibar_{|\Fbar}+ a \Ebar_{|\Fbar},
\end{equation}
where the divisors $\sMbar_{|\Fbar}$ and $\Ebar_{|\Fbar}$
are ample.
This is possible only if $\Fbar\iso \PP^2$, $\Xhat\iso\PP^1$,
$\Oh_{\Fbar}(\sMbar)\iso \Oh_{\Fbar}(\Ebar)\iso\Oh_{\PP^2}(1)$, and $a=1$.
From the exact sequence
\begin{equation}
0 \to \Oh_{\Xbar}(\sMbar-\Fbar)
\to \Oh_{\Xbar}(\sMbar)
\to \Oh_{\Fbar}(\sMbar)
\to 0
\end{equation}
we get
\begin{equation}
h^0(\Oh_{\Xbar}(\sMbar-\Fbar)) \ge
h^0(\Oh_{\Xbar}(\sMbar))
- h^0(\Oh_{\Fbar}(\sMbar))\ge2.
\end{equation}
Thus $\sMbar\ni \Fbar+\Lbar$, where $\Lbar\in |\sM-\Fbar|$ is a
mobile divisor. Hence there is a decomposition $-K_{X}\sim 2 F+2
L+\Xi$. In particular, $q\ge5$ and $F \qq L\qq A$. This implies that
$\fbar$ has no multiple fibers. So, the group $\Cl\Xbar$ is torsion
free. Since $\Oh_{\Fbar}(\Ebar)\iso\Oh_{\PP^2}(1)$, the class of
$\Ebar$ is not divisible in $\Cl\Xbar$. Hence $\Cl X$ is also torsion
free, a contradiction.

Therefore, $\sMbar$ is $\fbar$-vertical. Then $\sMbar=\fbar^*\sB$,
where $\sB$ is a linear system of Weil divisors on $\Xhat$ with $\dim
\sB\ge4$. We use the notation of~\ref{ss!Mfs}.
Let $\Gbar=\fbar^*\Theta$. We can write $\sB\qq t \Theta$ for some
$t\in \Z_{>0}$. Then
\begin{equation}
-K_{\Xbar}\qq 2 t\Gbar+\Xibar+a\Ebar,
\end{equation}
so $8\ge q\ge2t+1$ and $t\le3$. If $\Xhat\iso\PP^1$, we obviously have
$\dim\sB\le2$. Therefore, $\Xhat$ is a surface. Now we use
\ref{possibilities-Du-Val-del-Pezzo}.

If $t=1$, then $\dim\sB\le2$, a contradiction. Consider the case $t=2$.
Then $\dim\sB\ge4$ only in the case $\Xhat\iso\PP^2$. Then $q\ge5$,
$G\qq A$, and $m=2$. Since $\dim|G|\ge2$, by Lemma~\ref{Lemma-q=5} we
have $X\iso \PP(1^3,2)$. Consider the case $t=3$. Then $q\ge7$ and
$G\qq A$. Since $\dim \sB\ge4$, we have either $\Xhat\iso \PP^2$, $\Xhat\iso \PP(1,1,2)$, or $K_{\Xhat}^2=1$.
In either case $\dim|G|\ge1$ (recall that if $K_{\Xhat}^2=1$,
we take $\Theta=-K_{\Xhat}$).
By Lemma~\ref{Lemma-q=7} we get $X\iso \PP(1^2,2,3)$.

Now assume that $\fbar$ is birational.
We have
\begin{equation}
-K_{\Xhat}\sim 2 \sMhat+ \Xihat+ a \Ehat,
\end{equation}
where, as usual, we write $\Lahat=\fbar_*\Labar$ for the birational
transform of $\Xhat$ of a divisor (or a linear system) $\Labar$ on
$\Xbar$.

Clearly, $\dim\sMhat\ge\dim\sM$. If $(\Xhat,\sMhat)$ is not terminal, we
can repeat the procedure \ref{ss!Slink} and continue. Thus we may assume
that $(\Xhat,\sMhat)$ is as in \ref{p!dim-2M-q=13}--\ref{p!dim-2M-q=3}.
In particular, $\Cl\Xhat$ is torsion free and $\Xihat+ a \Ehat\sim
\Thhat$, where $\Thhat$ is the ample generator of $\Cl\Xhat$. So,
$\Xihat=0$, $a=1$, and $\Ehat\sim \Thhat$. In particular, the class of
$\wE$ is a primitive element of $\Cl\wX\iso \Z\oplus \Z$. In this case,
$\Cl X$ is also torsion free.
\end{proof}

\section{Proof of Main Theorem \ref{th!m}}\label{s!pf1}

Let $X$ be a $\Q$-Fano $3$-fold such that $-K_X\sim 2A$ for a primitive
element $A\in\Cl X$. Assume that $\dim|A|\ge4$ and  $K_X$ is not Cartier.
We apply Construction \ref{ss!Slink} with $\sM:=|A|$, $\la=2$ and $\Xi=0$.
By Lemma~\ref{lemma-del-pezzo} the pair $(X,\sM)$ is not terminal. Hence
in the notation of \eqref{eq!KX+laM}, the discrepancy $a>0$. On the
other hand, $a$ is an integer. Therefore, $a\ge1$.

\begin{lem}
The map $\fbar$ in \eqref{eq3.1} is birational.
\end{lem}

\begin{proof}
Suppose that $\fbar$ is not birational.
Let $\Fbar$ be a general fiber of $\fbar$.
If $\sMbar$ is $\fbar$-vertical, then $\sMbar=\fbar^*\sBhat$, where
$\sBhat$ is a linear system on $\Xhat$ whose class generates $\Cl\Xhat/\Tors$.
But then $\dim \sMbar=\dim \sBhat\le 2$ by \ref{possibilities-Du-Val-del-Pezzo},
contradicting our assumption.

Thus $\sMbar$ is $\fbar$-horizontal. Then
$-K_{\Fbar}=2\sMbar_{|\Fbar}+a\Ebar_{|\Fbar}$. This implies that
$\Fbar\iso \PP^2$, that is, $\fbar$ is generically a $\PP^2$-bundle
and $\Oh_{\Fbar}(\sMbar)\iso \Oh_{\PP^2}(1)$. From the exact
sequence
\begin{equation}
0 \to \Oh_{\Xbar}(\sMbar-\Fbar)
\to \Oh_{\Xbar}(\sMbar)
\to \Oh_{\Fbar}(\sMbar)
\to 0
\end{equation}
we get
\begin{equation}
h^0(\Oh_{\Xbar}(\sMbar-\Fbar)) \ge2.
\end{equation}
Therefore, $\sMbar\ni \Fbar+\Lbar$, where $\Fbar$ and
$\Lbar$ are mobile divisors. This contradicts $\qQ(X)=2$.
\end{proof}

\subsection{Notation}\label{notation-birational}
When $\fbar$ is birational, $\Xhat$ is a $\Q$-Fano. Recall that we write
$\Lahat=\fbar_*\Labar$ for the birational transform on $\Xhat$ of a
divisor (or a linear system) $\Labar$ on $\Xbar$. We have
\begin{equation}
\label{equation-q-main}
-K_{\Xhat}\sim 2 \sMhat+a\Ehat \quad\hbox{with} \quad a>0,\quad \dim \sMhat\ge4.
\end{equation}
By Proposition~\ref{p!dim-2M} the class of $\Ehat$ is the ample
generator of $\Cl\Xhat\iso \Z$, $\qhat=2n+1$, and $\sMhat\subset
|n\Ehat|$. Moreover, $\Xhat$ belongs to one of the possibilities listed
in Proposition~\ref{p!dim-2M}.

Assume first that $\qhat>3$. We consider the case $\qhat=3$ in the next
section. We make frequent use of the following easy observation.

\begin{rem}\label{Remark-decomposition}
In the notation of \ref{notation-birational}, assume that there is a
member $\Mhat\in \sMhat$ such that $\Mhat=\Lhat_1+\Lhat_2$, where
$\Lhat_1$ and $\Lhat_2$ are effective ample Weil divisors. Then either
$\Supp\Lhat_1= \Ehat$ or $\Supp\Lhat_2= \Ehat$.

Indeed, we can write
\begin{equation}
\sMbar\qq \Lbar_1+\Lbar_2 +\ga \Fbar,
\end{equation}
where $\Lbar_i$ is the birational transform of $\Lhat_i$ and $\ga\ge0$. Therefore,
\begin{equation}
\sM\sim f_*\chi^{-1}_*\sMbar\qq f_*\chi^{-1}_*\Lbar_1+f_*\chi^{-1}_*\Lbar_2+\ga F.
\end{equation}
Since the class of $A$ is a primitive element of $\Cl X$, we have either
$f_*\chi^{-1}_*\Lbar_1=0$ or $f_*\chi^{-1}_*\Lbar_2=0$ (and $\ga=0$).
\end{rem}

\begin{stheorem}{\bf Corollary.}\label{Corollary-decomposition}
Assume that we have $\dim|n\Ehat|=4$ in the notation of
\ref{notation-birational}. Then for any partition $n=n_1+n_2$, $n_i\in
\Z$ either $\dim|n_1\Ehat|\le 0$ or $\dim|n_2\Ehat|\le 0$.
\end{stheorem}

\begin{proof}
In this case $\sMhat=|n\Ehat|$ is a complete linear system.
Hence, one can take $\Lhat_i\in |n_i\Ehat|$.
\end{proof}

We consider the cases of Proposition~\ref{p!dim-2M} separately.

\begin{scase}{\bf Cases \ref{p!dim-2M-q=13},
\ref{p!dim-2M-q=9} and \ref{p!dim-2M-q=5a}.}
Then $\dim|n\Ehat|=4$ and $n$ is even. Apply Corollary \ref{Corollary-decomposition}
with $n_1=n_2=n/2$.
We get a contradiction because $\dim|n_i\Ehat|>0$.
\end{scase}

\begin{scase}{\bf Case \ref{p!dim-2M-q=11}, $\Xhat\iso \PP(1,2,3,5)$.}
Then $n=5$ and $\dim|n\Ehat|=5$. Thus $\sMhat\subset|5\Ehat|$ is a
subsystem of codimension $\le 1$. Since $\dim 2\Ehat|=1$, we can take
$\Lhat_1 \in |2\Ehat|$ so that $\Lhat_1 \ne2\Ehat$. Since
$\dim|3\Ehat|=2$, there exists a one dimensional family of divisors
$\Lhat_2 \in |3\Ehat|$ such that $\Lhat_1+\Lhat_2\in \sMhat$. So we may
assume that $\Lhat_2 \ne3\Ehat$. By Remark \ref{Remark-decomposition} we
get a contradiction.
\end{scase}

\begin{scase}{\bf Case \ref{p!dim-2M-q=7}}, that is, $\Xhat \iso
\PP(1^2,2,3)$. Then $n=3$ and $\dim|n\Ehat|=6$. Thus $\sMhat\subset
|3\Ehat|$ is a subsystem of codimension $\le 2$. Since $\dim|\Ehat|=1$,
|we can take $\Lhat_1 \in|\Ehat|$ so that $\Lhat_1 \ne\Ehat$. Since
$\dim|2\Ehat| = 3$, there exists a one dimensional family of divisors
$\Lhat_2 \in |2\Ehat|$ such that $\Lhat_1 + \Lhat_2 \in \sMhat$. So we
may assume that $\Lhat_2 \ne2\Ehat$. By Remark
\ref{Remark-decomposition} we get a contradiction. \end{scase}

\begin{scase}{\bf Case \ref{p!dim-2M-q=5}, $\Xhat\iso
\PP(1^3,2)$.} Then $n=2$ and $\dim|n\Ehat|=6$. Thus $\sMhat\subset
|2\Ehat|$ is a subsystem of codimension $\le 2$.

Assume that $\fbar(\Fbar)$ is a curve. Then
\begin{equation*}
K_{\Xbar}=\fbar^*K_{\Xhat}+ \Fbar\quad\hbox{and}\quad
\Ebar=\fbar^*\Ehat- \ga\Fbar.
\end{equation*}
Since any member of $|\Ehat|$ is smooth in codimension one, $\ga\le 1$.
Moreover, since $n\Ebar$ is not mobile for any $n$, we have $\ga>0$. Hence,
$\ga=1$. So,
\begin{equation*}
K_{\Xbar}+5\Ebar+ 4\Fbar=\fbar^*(K_{\Xhat}+5\Ehat)\sim 0.
\end{equation*}
This implies that $-K_X$ is divisible by $4$, a contradiction.

Hence $\fbar(\Fbar)\in \Xhat$ is a point, say $\Phat$. If $\Phat\in
\Xhat$ is the point of index $2$, then $\fbar$ is the blowup of the
maximal ideal \cite{Kawamata-1996}. In this case $\Xbar$ has exactly two
extremal contractions: $\fbar$ and the $\PP^1$-bundle induced by the
projection $\PP(1^3,2)\broken \PP^2$. On the other hand, the second
contraction must be birational, a contradiction. Hence $\Pbar\in \Xhat$
is a smooth point.

Let $\sLhat:=|\Ehat|$. Take a general member $\Lhat_1\in \sLhat$.
Dimension counting shows that there exists $\Lhat_2\in \sLhat$ such that
$\Lhat_1+\Lhat_2\in \sMhat$. If $\Lhat_2\ne\Ehat$, we get a
contradiction by Remark \ref{Remark-decomposition}. Thus $\Lhat_2=\Ehat$
for any choice of $\Lhat_1\in \sLhat$. Therefore, $\Ehat+\sLhat\subset
\sMhat$ and we can write $\sMbar\qq \sLbar+\Ebar+\ga\Fbar$, where
$\ga\ge0$. Then
\begin{equation}
0\sim K_{\Xbar}+2\sMbar+\Ebar\qq K_{\Xbar}+2\sLbar+3\Ebar+2\ga\Fbar.
\end{equation}
Note that the only base point of $\sLhat$ is the point of index $2$.
Hence, $\sLbar\qq \fbar^*\sLhat$. Let $\sLhat'\subset \sLhat$ be the
subsystem consisting of elements passing through $\Phat$. Then we can
write
\begin{equation}
\sLbar'\qq\fbar^*\sLhat'-\de\Fbar\qq \sLbar-\de\Fbar,
\quad\hbox{with $\de>0$.}
\end{equation}
Therefore,
\begin{equation}
0 \qq K_{\Xbar}+2\sLbar+3\Ebar+2\ga\Fbar\qq
K_{\Xbar}+2\sLbar'+3\Ebar+2(\de+\ga)\Fbar.
\end{equation}
This gives us $-K_X\qq 2\sL'+2(\de+\ga)F$ which contradicts $\qQ(X)=2$.
\end{scase}

\section{Conclusion of the proof of Main Theorem \ref{th!m}}\label{s!pf2}
This section considers Case~\ref{p!dim-2M-q=3}, when $\Xhat=Q\subset
\PP^4$ is a smooth quadric. Then $\sMhat=|\Oh_Q(1)|$ is a complete
linear system, and in particular is base point free. Thus $\sMbar\qq
\fbar^*\sMhat$. We also have $\Ehat\in|\Oh_Q(1)|$ and
$\fbar(\Fbar)\subset\Ehat$.

\begin{lem}
$\Ga:=\fbar(\Fbar)$ is a curve.
\end{lem}
\begin{proof}
Assume that $\fbar(\Fbar)$ is a point. Let $\sMhat'\subset \sMhat$
be the sub\-system consisting of elements passing through $\fbar(\Fbar)$.
Then we can write
\begin{equation}
\sMbar' \qq \fbar^*\sMhat'-\de\Fbar\qq \sMbar-\de\Fbar,
\quad\hbox{with $\de>0$.}
\end{equation}
Therefore,
\begin{equation}
\begin{aligned}
0 &\qq \fbar^*(K_{\Xhat}+2\sMhat'+\Ehat) \qq \fbar^*(K_{\Xhat}+2\sMhat+\Ehat) \\
&\qq K_{\Xbar}+2\sMbar+\Ebar\qq K_{\Xbar}+2\sMbar'+\Ebar+2\de\Fbar.
\end{aligned}
\end{equation}
This gives us $-K_X\qq\, 2\sM'+2\de F$ which contradicts $\qQ(X)=2$.
\end{proof}

\begin{lem}\label{lem!Ebar}
$\Ebar\iso \PP(1,1,2)$, $\FF_2$, or\/ $\PP^1\times\PP^1$.
\end{lem}
\begin{proof}
Clearly, $\Ehat\iso \PP(1,1,2)$ or $\PP^1\times \PP^1$. In particular,
the pair $(\Xhat,\Ehat)$ is plt. Since $K_{\Xbar}\qq
\fbar^*K_{\Xhat}+\Fbar$ and $\Ehat$ is smooth at the generic point of
$\Ga$, we have
\begin{equation}
\label{eq!q=3-KE}
K_{\Xbar}+\Ebar\sim \fbar^*(K_{\Xhat}+\Ehat).
\end{equation}
Hence the pair $(\Xbar,\Ebar)$ is plt and the divisor $K_{\Xbar}+\Ebar$
is Cartier. By adjunction, the surface $\Ebar$ has at worst Du Val
singularities. Moreover, $K_{\Ebar}=\fbar_{|\Ebar}^*K_{\Ehat}$, that is,
the restriction $\fbar_{\Ebar}$ is either an isomorphism or the minimal
resolution of $\Ehat$.
\end{proof}
\begin{lem}
$-K_{\Xbar}$ is nef.
\end{lem}

\begin{proof}
Recall that by our construction $\Xbar$ has exactly two extremal rays.
Denote them by $R_1$ and $R_2$. One of them, say $R_1$, is generated by
nontrivial fibers of $\fbar$. Let $C$ be an extremal curve on $\Xbar$
that generates $R_2$. Assume that $-K_{\Xbar}$ is not nef. Then
$K_{\Xbar}\cdot C>0$ and $C$ must be a flipped curve (because a
divisorial contraction must be $K$-negative in our situation). Since
$-K_{\Xbar}\qq \Ebar+2\fbar^*\Ehat$, we have $\Ebar\cdot C<0$. In
particular, $C\subset \Ebar$. Since $C$ is a flipped curve, it cannot be
mobile on $\Ebar$, that is, $\dim|C|=0$. By Lemma~\ref{lem!Ebar}, the
only possibility is that $\Ebar\iso \FF_2$ and $C$ is the negative
section of $\FF_2$. But in this case $C$ is contracted by $\fbar$ to a
point, that is, the class of $C$ lies in $R_1$, a contradiction.
\end{proof}

\begin{lem} $K_{\Xbar}$ is not Cartier at some point of\/ $\Ebar$.
\end{lem}
\begin{proof}
By \eqref{eq!q=3-KE} the divisor $K_{\Xbar}$ is Cartier outside $\Ebar$.
Assume that $K_{\Xbar}$ is Cartier near $\Ebar$. Since $-K_{\Xbar}$ is
nef, the map $\Xbar\broken \wX$ is either an isomorphism or a flop. In
either case $\wX$ has the same type of singularities as $\Xbar$, that is,
$K_{\wX}$ is Cartier. By the classification of extremal contractions of
Gorenstein terminal $3$-folds \cite{Cutkosky-1988} the divisor $2K_X$ is
Cartier. This contradicts the following remark.
\end{proof}

\begin{stheorem}{\bf Corollary.}
The curve $\Ga$ has a singular point that is not a local complete
intersection.
\end{stheorem}

\begin{proof}
Indeed otherwise by \cite[Prop.\ 4.10.1]{KM92} the map $\fbar$ is the
blowup of $\Ga$ and $K_{\Xbar}$ is Cartier.
\end{proof}

\begin{stheorem}{\bf Corollary.}
$\Ehat\iso \PP(1,1,2)$, the curve $\Ga$ is not a Cartier divisor on
$\Ehat$, and $\Ga$ is singular at the vertex of $\PP(1,1,2)$.
\end{stheorem}

\begin{lem}
$\deg \Ga=5$.
\end{lem}
\begin{proof}
Let $\Chat\subset \Ehat$ be a general hyperplane section. Since
$-K_{\Xbar}$ is nef,
\begin{equation}
0\le -K_{\Xbar}\cdot \Cbar= -K_{\Xhat}\cdot \Chat - (\Ga \cdot
\Chat)_{\Ehat}=6-\deg \Ga.
\end{equation}
Since $\Ga$ is not a Cartier divisor on $\Ehat$, its degree should be
odd. If $\deg\Ga\ne5$, then $\Ga$ is either a line or a twisted cubic.
In particular, it is smooth, a contradiction.
\end{proof}

\begin{case}
Thus $\deg\Ga=5$ and $\Ga$ is singular. Then $\Ga$ can be given, in
coordinates $u_1,u_2,v$ for $E\iso \PP(1,1,2)$, by an equation
$\ga=v\al_3+\be_5$, where $\al_3(u_1,u_2)$ and $\be_5(u_1,u_2)$ are
homogeneous polynomial of the indicated degree. Thus $P$ is a triple
point of $\Ga$ and is its only singularity. Thus $\Ga$ is as in Main
Example~\ref{ss!MEx} and the conclusion of Theorem~\ref{th!m}. By
\cite[Th.\ 4.9]{KM92} the extraction $\fbar\colon \Xbar\to Q=\Xhat$ of
$\Ga$ is unique up to isomorphism over $Q$. Since $\rho(\Xbar/Q)=1$, the
Sarkisov link \eqref{eq!Sark} is uniquely determined. This completes the
proof of Theorem~\ref{th!m}.
\end{case}

\section{Examples}\label{s!mEx}

\subsection{Symbolic blowup}\label{ss!sBl}

This section is closely related to parts of Tom Ducat's thesis
\cite{Du15}, and we acknowledge his help with our treatment.

Let $\Ga\subset M$ be a reduced singular curve in a nonsingular $3$-fold.
The {\em symbolic blowup} of $\Ga$ in $M$ is the relative Proj of the {\em
symbolic algebra}, the graded algebra $\bigoplus_{n\ge0} \sI_\Ga^{[n]}$,
where $\sI_\Ga^{[n]}$ is the $n$th symbolic power, that is, the ideal in
$\Oh_Q$ of functions vanishing $n$ times at the generic point of $\Ga$. In
other words, in the primary decomposition of $\sI_\Ga^n$, ignore the
embedded component at singular points $P$ of $\Ga$. (Primary decomposition
is built into the computer algebra packages.)

In our case, $\Ga$ is a curve contained in a $\frac12(1,1)$ orbi\-fold
point $P\in E_0\subset M$ as a Weil divisor, whose class generates the
local class group $\Cl_PE_0\iso\Z/2$. For simplicity, we treat
$\Ga\subset E_0\subset M$ as germs around a singular point $P\in\Ga$ in
local analytic coordinates (but see \ref{ss!alt}). Write
$\C^2_{\Span{u_1,u_2}}$ for the orbifold double cover of $P\in E_0$, and
\begin{equation}
(x_1,x_2,x_3)=(u_1^2,u_1u_2,u_2^2)
\end{equation}
for the invariant monomials. Then $M$ has local coordinates
$x_1,x_2,x_3$, with $g=x_1x_3-x_2^2$ the local equation of $E_0$, and
$\Ga\subset E_0$ corresponds to an invariant curve
$\Ga:(\ga=0)\subset\C^2$, with equation $\ga=\ga(u_1,u_2)$ an odd
function of the orbinates.

To see $\Ga\subset M$ in equations, first render into $x_i$ the invariant
multiples $u_1\ga$ and $u_2\ga$ of $\ga$, say as
\begin{equation}
u_1\ga=bx_1-ax_2=-f_2 \quad\hbox{and}\quad u_2\ga=bx_2-ax_3=f_1,
\end{equation}
where $a,b$ are functions of $x_1,x_2,x_3$. Taking into account that $a$
and $b$ are in the maximal ideal $(x_1,x_2,x_3)$ (because the curve
$\Ga$ is singular at $P$, and not locally planar), and with a little
massaging, we can put the generators of $\sI_\Ga$ and the syzygies
between them in the determinantal form $\bigwedge^2M=(f_1,f_2,g)$, where
\begin{equation}
M = \begin{pmatrix}
x_1 & x_2 & a_2x_2+a_3x_3 \\
x_2 & x_3 & b_1x_1+b_2x_2
\end{pmatrix}
\quad\hbox{and}\quad
M\begin{pmatrix} f_1\\ f_2\\ g \end{pmatrix}\equiv 0.
\label{eq!Mf=0}
\end{equation}

In this case, the symbolic algebra needs just one further generator in
degree~2, whose restriction to $E_0$ is the local equation
\begin{equation}
b^2x_1-2abx_2+a^2x_3
\label{eq!2C}
\end{equation}
of the Cartier divisor $2\Ga\subset E_0$. Rather than primary
decomposition, we derive this final generator and its relations by
unprojection.

Replacing $f_1,f_2,g$ by $\xi_1,\xi_2,\eta$ in \eqref{eq!Mf=0} gives
\begin{equation}
\begin{pmatrix}
x_1 & x_2 & a_2x_2+a_3x_3 \\
x_2 & x_3 & b_1x_1+b_2x_2
\end{pmatrix}
\begin{pmatrix} \xi_1\\ \xi_2\\ \eta \end{pmatrix}=0.
\label{eq!Mxi}
\end{equation}
Equations \eqref{eq!Mxi} define the blowup of the ideal
$\sI_\Ga=(f_1,f_2,g)$ as a codimension~2 complete intersection in
$M\times\PP^2_{\Span{\xi_1,\xi_2,\eta}}$, containing the ``irrelevant''
codimension~3 complete intersection $V(\xi_1,\xi_2,\eta)$. However, $M$
has entries in $(x_1,x_2,x_3)$, so it also contains the codimension~2
complete intersection $V(x_1,x_2,x_3)$ -- the blowup of $\Ga$ must
contain $\PP^2$ over the origin (because $\Ga$ is not a local complete
intersection). We rearrange
\eqref{eq!Mxi} as
\begin{equation}
\begin{pmatrix}
\xi_1 & \xi_2+a_2\eta & a_3\eta \\
b_1\eta & \xi_1+b_2\eta & \xi_2
\end{pmatrix}
\begin{pmatrix}
x_1 \\ x_2 \\ x_3
\end{pmatrix}=0.
\end{equation}
The unprojection of $V(x_1,x_2,x_3)$ is given by the $4\times4$
Pfaffians of
\begin{equation}
\begin{pmatrix}
\ze & \xi_1 & \xi_2+a_2\eta & a_3\eta \\
& b_1\eta & \xi_1+b_2\eta & \xi_2 \\
&& x_3 & -x_2 \\
&&& x_1
\end{pmatrix}.
\label{eq!ze}
\end{equation}

Geometrically, this is the blowup of $\sI_\Ga$ followed by the
unprojection contracting $\PP^2$ over the origin to a point.
We can also view it simply as a practical means of writing
down the generator $h$ in degree 2 satisfying
\begin{equation}
(x_1,x_2,x_3)h = \bigwedge\nolimits^2 \begin{pmatrix}
f_1 & f_2+a_2g & a_3g \\
b_1g & f_1+b_2g & f_2
\end{pmatrix},
\end{equation}
so that clearly $h\in\sI_\Ga^{[2]}$, without computer algebra. In
computational terms, this means that we can modify $f_1^2, f_1f_2,
f_2^2$ modulo $g\cdot\sI_\Ga$ to make them identically divisible by
$x_3,x_2,x_1$ respectively, with (say)
\begin{equation}
f_1f_2-b_1a_3g^2 = -x_2h,
\end{equation}
where $h_{|E_0}$ is the equation \eqref{eq!2C} defining $2\Ga\subset E_0$.

\begin{prop} The symbolic algebra of\/ $\sI_\Ga$ is generated by
$\xi_1,\xi_2$, $\eta$ in degree~$1$ (corresponding to $f_1,f_2,g$), and
$\ze$ in degree~$2$ (corresponding to $h$). The ideal of relations is
generated by the maximal Pfaffians of \eqref{eq!ze}.

Thus the symbolic blowup $M_1\to M$ of $\Ga$ is the codimension~$3$
Gorenstein subvariety
\begin{equation}
M_1\subset M\times\PP(1,1,1,2)_{\Span{\xi_1,\xi_2,\eta,\ze}}
\end{equation}
defined by the Pfaffians of \eqref{eq!ze}. It has the following
properties. If $\Ga$ is nonsingular at $P$ it is not applicable. If $\Ga$ is
singular it defines a morphism $M_1\to M$ which is the ordinary blowup
of $\sI_\Ga$ outside the origin. The birational transform $E_1\subset M_1$
is isomorphic to $E_0$.

The fibre of $M_1$ over the origin is $\PP(2,1)_{\Span{\ze,\eta}}$
passing through $P_\ze$, which is a $\frac12(1,1,1)$ orbifold point, and
at most one more singular point. If $\mult_P\Ga\ge5$ then $P_\eta\in
M_1$ has embedding dimension~$3$, so is not terminal. If $\mult_P\Ga=3$
then $M_1$ is terminal, and in fact:
\begin{enumerate}
\renewcommand{\labelenumi}{(\arabic{enumi})}
\renewcommand\theenumi{(\arabic{enumi})}

\item $M_1$ is quasismooth if $\Ga$ has $3$ distinct tangent branches.

\item $M_1$ has a c$A_1$ point if $\Ga$ has a double tangent branch.

\item $M_1$ has a c$A_2$ point if $\Ga$ has a triple tangent branch.

\end{enumerate}
\end{prop}

In the local description, the c$A_1$ and c$A_2$ points of $M_1$ are arbitrary.

\medskip
\paragraph{\em Proof.\ }
Although the precise statement is somewhat involved, the proof is easy.
The generators and relations follow from the hyperplane section
principle: indeed, the symbolic algebra restricted to $E_0$ is just an
algebra of $\Z/2$ invariants, generated by $u_1\ga,u_2\ga$ and $\ga^2$,
and the restriction map is surjective by our choice of generators
$\xi_1,\xi_2$ and $\ze$.

The birational transform $E_1\to E_0$ is an isomorphism because the
symbolic algebra of the singular (nonplanar) curve $\Ga\subset M$ maps
onto that of the $\Q$-Cartier divisor on $\Ga\subset E_0$. The analysis
of the singularities is straightforward. The cases correspond to the
different possibilities for the cubic leading terms in $\ga$ coming from
\begin{equation}
\hphantom{\qed}\enspace
b_1u_1^3+b_2u_1^2u_2-a_2u_1u_2^2-a_3u_2^3.
\qquad\qed
\end{equation}

\begin{rem}
In the case that $\Ga$ has distinct tangent branches, its symbolic blowup
can be done as an explicit construction in the nonsingular category that
is folklore in the subject: blowing up $P$ gives an exceptional
$\Pi=\PP^2$ with normal bundle $\Oh(-1)$. The 3 branches of $\Ga$ meet
$\Pi$ at noncollinear points, and blowing up $\Ga$ produces a dP$_6$ with
a hexagon formed of the 3 blown up lines on $\Pi$ together with the 3
lines joining them in $\Pi$, that are $(-1,-1)$ curves. Flopping these
takes $\Pi$ into $\Pi'=\PP^2$ with normal bundle $\Oh(-2)$ by a standard
quadratic transformation $\Pi\broken\Pi'$, and $\Pi'$ contracts to a
$\frac12(1,1,1)$ point.
\end{rem}

\subsection{Two examples} \label{ss!Ex2}
We apply this to contruct two families of $\Q$-Fano $3$-folds $X$ and
$Y$ of index~$2$ with $\Cl=\Z\cdot A$, each with a single $\frac13(1,2,2)$
orbifold point and invariants
\begin{equation}
\begin{array}{llll}
-K_X=2A_X & \hbox{with} & A_X^3=\tfrac{10}{3}, & \dim|A_X|=4, \\[4pt]
-K_Y=2A_Y & \hbox{with} & A_Y^3=\tfrac{7}{3}, & \dim|A_Y|=3.
\end{array}
\end{equation}
Their Hilbert series come from this by the Ice Cream formula
of \cite{Ice}:
\begin{equation}
\renewcommand{\arraycolsep}{.2em}
\begin{array}{rcccl}
P_{X,A_X}(t) &=& \tfrac{1+t+t^2}{(1-t)^4} + \tfrac{t^2}{(1-t)^3(1-t^3)}
 &=& \tfrac{1+2t+4t^2+2t^3+t^4}{(1-t)^3(1-t^3)} \\[12pt]
P_{Y,A_X}(t) &=& \tfrac{1+t^2}{(1-t)^4} + \tfrac{t^2}{(1-t)^3(1-t^3)}
 &=& \tfrac{1+t+3t^2+t^3+t^4}{(1-t)^3(1-t^3)}.
\end{array}
\label{eq!PY}
\end{equation}

\begin{exa} \label{ex!ex2}
Let $E_0\subset\PP^3$ be the ordinary quadratic cone and $\Ga_7\subset
E_0\subset\PP^3$ a curve of degree~7 that is singular at the node $P\in
E_0$, with $\mult_P\Ga_7=3$. The symbolic blowup of $\Ga_7$ defines an
extremal extraction $Y_1\to\PP^3$, with the birational transform
$E_1\subset Y_1$ isomorphic to $E_0$ and to $\PP(2,1,1)$.

Write $B$ for the polarizing $\Oh(1)$ of $\PP^3$ and its pullback to
$Y_1$, so that $E_0\sim 2B$, and let $F\subset Y_1$ be the scroll over
$\Ga_7$. Note that $2\Ga_7\sim 7B$ in $\Pic E_0$, so $2F\sim7B$ in $\Pic
E_1$. In $\Cl Y_1$ we have
\begin{equation}
2B=E_1+F \quad\hbox{and}\quad K_{Y_1}=-4B+F=-2B-E_1.
\label{eq!G1}
\end{equation}

We give $Y_1$ the polarizing divisor $A_1=B+\frac23E_1$. Then
\begin{equation}
3A_1=3B+2E_1=7B-2F
\end{equation}
so that $3A_1$ is a Cartier divisor restricting to a linearly trivial
divisor on $E_1$, that is, $\Oh_{E_1}(3A_1)\iso\Oh_{E_1}$. Standard use
of vanishing gives that $H^0(Y_1,\Oh_{Y_1}(3A_1))\onto H^0(\Oh_{E_1})$
is surjective, so that $|3A_1|$ is a free linear system, ample outside
$E_1$, so contracts $E_1$ to a $\frac13(2,1,1)$ orbifold point on a
$3$-fold $Y$.

Also, \eqref{eq!G1} gives
\begin{equation}
K_{Y_1}-\tfrac13E_1=-2B-\tfrac43E_1=-2A_1.
\end{equation}
Hence $-K_Y=2A$, with $A$ an ample Weil divisor on $Y$. The contraction
$Y_1\to Y$ is the Kawamata blowup of $\frac13(2,1,1)$, with discrepancy
$\tfrac13E_1$.
\end{exa}

\begin{exa} \label{ex!main}
The Main Example of Theorem~\ref{th!m} is almost the same. We start
from the nonsingular quadric $Q\subset\PP^4$ and the ordinary quadratic
cone obtained as the intersection $E_0=T_{P,Q}\cap Q$ with its tangent
hyperplane at a point $P\in Q$. Let $\Ga_5\subset E_0$ be a irreducible
quintic curve, assumed singular at $P$ (it follows that $\mult_P\Ga=3$).

The symbolic blowup $X_1\to Q$ of $\Ga_5$ has exceptional scroll $F$, and
birational transform $E_1\iso E_0$. As before, write $B=\Oh(1)$ for the
polarizing divisor of $Q$, so that $E_0\sim B$, and also for its
pullback to $X_1$. Thus in $\Cl X_1$ we have
\begin{equation}
B=E_1+F \quad\hbox{and}\quad K_{X_1} = -3B+F = -2B-E_1.
\end{equation}

We give $X_1$ the polarising divisor $A_1 = B+\tfrac23 E_1$.
Then $3A_1=3B+2E_1=5B-2F$ is a Cartier divisor with
$\Oh_{E_1}(3A_1)\iso\Oh_{E_1}$ with surjective restriction
$H^0(\Oh_{X_1}(3A_1))\onto H^0(\Oh_{E_1})$. Thus $|3A_1|$ is a
free linear system contracting $E_1$ to a $\frac13(1,2,2)$
point. Now $K_{X_1}=-2B-G1=-2A_1+\frac13E_1$ so that $-K_X=2A$
with $A$ an ample Weil divisor, and $X_1\to X$ is the Kawamata
blowup, with discrepancy $\frac13E_1$.
\end{exa}

\subsection{Alternative graded ring constructions} \label{ss!alt}
We can treat the examples of \ref{ss!Ex2} in graded ring terms. This is
how we originally discovered them. Moreover, the algebra is interesting
in its own right, and displays features that are possibly typical for
index~2 Fano constructions.

The construction of $Y$ is immediate. Its Hilbert series \eqref{eq!PY} is
\begin{equation}
\frac{1-2t^2-3t^3+3t^4+2t^5-t^7}{(1-t)^4(1-t^2)^2(1-t^3)},
\end{equation}
indicating the codimension~3 subvariety
$Y\subset\PP(1^4,2^2,3)_{\Span{x_{0\dots3},y_1,y_2,z}}$ defined by the
maximal Pfaffians of a $5\times5$ matrix of degrees
\begin{equation}
\left(
\begin{smallmatrix}
3&2&2&2 \\
&2&2&2 \\
&&1&1 \\
&&&1
\end{smallmatrix}\right),
\quad\hbox{typically}\quad
\begin{pmatrix}
z & y_1 &y_2+a_2& a_3 \\
& b_1 &y_1+b_2& y_2 \\
&& x_3&-x_2 \\
&&&x_1
\end{pmatrix},
\label{eq!Pfz}
\end{equation}
with $a_2,a_3,b_1,b_2$ general quadratic forms in $x_0,x_1,x_2,x_3$. Every
$Y$ in this family is given in this way.

One can follow the argument back to see that in this case
$\Ga_7\subset\PP^3$ is defined by $\bigwedge^2M=0$ where
\begin{equation}
M = \begin{pmatrix}
x_1 & x_2 & a_2x_2+a_3x_3 \\
x_2 & x_3 & b_1x_1+b_2x_2
\end{pmatrix},
\end{equation}
and the $y_i$ in degree~2 are the rational forms solving
\begin{equation}
M \begin{pmatrix}
y_1 \\ y_2 \\ 1
\end{pmatrix} = 0.
\end{equation}

Plausible though it may seem at first sight, it is a mistake to confuse
the codimension~2 variety $\Ybar_{4,4}\subset\PP(1^4,2^2)$ defined by
these equations with the symbolic blowup $Y_1$ of $\Ga_7$. The latter is a
relative construction over $\PP^3$, and is contained in
$\PP^3\times\PP^2$, so it has the ratios $f_1:f_2:g$ as regular
functions, where $g=x_1x_3-x_2^2$ is the equation of $E_0$. It is not
simply polarized or projectively Gorenstein. As we have seen, $Y_1$ has
just one orbifold point of type $\frac12(1,1,1)$.

In contrast, $\Ybar_{4,4}$ contains $\PP(1,2,2)_{\Span{x_0,y_1,y_2}}$
with ideal $(x_1,x_2,x_3)$, and has $\PP^1_{\Span{y_1,y_2}}$ as
$\frac12(1,1)$ orbifold locus, so is not terminal. It is clearly
obtained from $Y\subset\PP(1^4,2^2,3)$ by eliminating $z$. Putting back
$z$ is a Type~I or Kustin--Miller unprojection, with the Pfaffians of
\eqref{eq!Pfz} giving the linear relations for $z$, so is perfectly
valid as a construction of $Y$. However, the birational relation between
$Y$ and $\Ybar_{4,4}$ involves first the weighted blowup of the
$\frac13(1,2,2)$ point with the given weights $(1,2,2)$, not the
Kawamata blowup with weights $(2,1,1)$, and this takes us outside the
Mori category. A similar thing happens in many other constructions or
attempted constructions of index~2 $\Q$-Fanos.

There is a similar narrative for the Main Example $X$, starting from
$\Ga_5\subset E_0\subset Q$. The Hilbert series $P_X$ has the form
$\frac{N(t)}{(1-t)^5(1-t^2)^2(1-t^3)}$ with numerator
\begin{equation}
\begin{aligned}
N(t)=1-t^2-4t^3&-4t^4 \\
&+4t^4+8t^5+4t^6 \\
&\hphantom{+4t^4+8t^5\ }-4t^6-4t^7-t^8+t^{10}.\\
\end{aligned}
\end{equation}
We keep the masked terms $-4t^4+4t^4$ to indicates that $R(X,A)$ needs
4~relations in degree~4. In fact, in order to have
$\frac13(1,2,2)_{x_4,y_1,y_2}$ at $P_z$, there must be 4~relations
$zx_i=c_i$ to eliminate $x_0,\dots,x_3$ there.

Now eliminating $z$ projects $X$ to $\Xbar\subset\PP(1^5,2^2)$ in
codimension~3, with the Hilbert series
\begin{equation}
\frac{1-t^2-4t^3+4t^4+t^5-t^7}{(1-t)^5(1-t^2)^2},
\end{equation}
which corresponds to the Pfaffians of a skew $5\times5$ matrix of
degrees
\begin{equation}
\left(
\begin{smallmatrix}
1&1&1&2 \\
&1&1&2 \\
&&1&2 \\
&&&2
\end{smallmatrix}\right),
\quad\hbox{typically}\quad
\begin{pmatrix}
x_0 & x_1 & x_2 & a \\
& x_2 & x_3 & b \\
&& x_4& -y_1 \\
&&& y_2
\end{pmatrix}.
\label{eq!Pfz2}
\end{equation}
Here we choose coordinates on $\PP^4$ with $Q:x_0x_4-x_1x_3+x_2^2=0$ and
$P=(1,0,\dots,0)$, making $T_{P,Q}:x_4=0$ and $E_0:x_1x_3=x_2^2$. Let
$\Ga_5\subset E_0\subset\PP^3$ be an irreducible curve of degree~5,
assumed to be singular at $P$.

In \eqref{eq!Pfz2}, $a$ and $b$ are quadratic forms in $x_{0\dots4}$.
The conditions that $\Xbar$ defined by the Pfaffians of \eqref{eq!Pfz2}
contains $\PP(1,2,2)_{\Span{x_4,y_1,y_2}}$ defined by the ideal
$(x_0,\dots,x_3)$ is that $a,b$ do not contain $x_4^2$. Then the 7
entries in the first two rows of \eqref{eq!Pfz2} are in the ideal
$(x_0,x_1,x_2,x_3)$, so it is a Jerry$_{12}$. At the same time, the
equations of $\Ga_5\subset E_0$ take the determinantal format
$\bigwedge^2M$ with $M$ as in \eqref{eq!Mf=0}.

As before, unprojecting $\PP(1,2,2)_{\Span{x_4,y_1,y_2}}\subset\Xbar$ is
a contruction of $X$ as a Type~I unprojection from $\Xbar\supset
\PP(1,2,2)$, but $\Xbar$ itself again has a line of $\frac12$ points, so
is not Mori category.

The ``double Jerry'' calculations of \cite[9.2]{TandJ} gives the
unprojection variable $z$ and most of its unprojection equations
$zx_i=c_i$. Eliminating the pivot $m_{12}=x_0$ from the Pfaffians of
\eqref{eq!Pfz2}, gives two equations
\begin{equation}
\begin{pmatrix}
x_1 & x_2 & a \\
x_2 & x_3 & b
\end{pmatrix}
\begin{pmatrix}
y_2 \\
y_1 \\
x_4
\end{pmatrix}=0
\quad\hbox{with}\quad
\begin{array}{l}
a=a_2x_2+a_3x_3,\\ b=b_1x_1+b_2x_2,
\end{array}
\end{equation}
that we rearrange
\begin{equation}
\begin{pmatrix}
b_1x_4 & y_2+b_2x_4 & y_1 \\
y_2 & y_1+a_2x_4 & a_3x_4
\end{pmatrix}
\begin{pmatrix}
x_1 \\
x_2 \\
x_3
\end{pmatrix}=0.
\label{eq!2J}
\end{equation}
From this we assemble a second Jerry$_{12}$ matrix 
\begin{equation}
\begin{pmatrix}
-z & b_1x_4 & y_2+b_2x_4 & y_1 \\
& y_2 & y_1+a_2x_4 & a_3x_4 \\
&& x_3& -x_2 \\
&&& x_1
\end{pmatrix},
\label{eq!aPf}
\end{equation}
whose maximal Pfaffians provide the equations for $x_1z,x_2z,x_3z$.
The final unprojection equation for $x_0z$
\begin{equation}
-x_0z = (b_1x_1 + b_2x_2 -a_2x_3)y_1 + a_3(x_3y_2 + b_1x_2x_4 + b_2x_3x_4).
\end{equation}
exists by the theory of Kustin--Miller unprojection, but we don't know
any smart way of deducing it. It has to be calculated by a laborious
primary decomposition or colon ideal calculation, or by writing out the
Kustin--Miller complexes.

The Jerry$_{12}$ matrix \eqref{eq!aPf} defines a codimension~3
subvariety in the family of our second example $Y\subset\PP(1^4,2^2,3)$
(compare \eqref{eq!Pfz}), but specialized to contain
$\PP^2_{\Span{x_1,x_2,x_3}}$ defined by the codimension~4 ideal
$(x_4,y_1,y_2,z)$. This is a third construction of our Main Example $X$.

To do this from scratch: in Example~\ref{ex!ex2}, suppose that the curve
$\Ga_7\subset E_0\subset\PP^3$ break up as the plane conic section
$(x_0=0)$ plus a quintic $\Ga_5$. Blowing up $\Ga_7$ transforms the
plane $(x_0=0)$ into a copy of $\PP^2$ with normal bundle $\Oh(-1)$,
that contracts to a point of the quadric $Q\subset\PP^4$, taking $E_0$
and $\Ga_5\subset E_0$ isomorphically into the data for \ref{ex!main}.

\subsection{Summary: Three constructions of $X$}
Our Main Example $X$ can be obtained in three different ways

\begin{enumerate}
\item The symbolic blowup of $\Ga_5\subset E_0\subset Q$ followed by the
contraction of $E_1$. Viewed from $X$, this is the Sarkisov link from
its $\frac13(1,2,2)$ point $P$ of Section~\ref{ss!MEx}; it is initiated
by the Kawamata blowup, that is the $(2,1,1)$ weighted blowup of $P$.

\item Construct the codimension~3 variety $\Xbar\subset\PP(1^5,2^2)$
given in the Pfaffian form \eqref{eq!Pfz2}, containing $\PP(1,2,2)$,
then unproject this plane. Viewed from $X$, this starts from
the $(1,2,2)$ weighted blowup of $P$, which introduces a line of
$\frac12(1,1)$ orbifold points, so takes us out of the Mori category.

\item Construct the codimension~3 variety $Y'\subset\PP(1^4,2^2,3)$ as
in Example~\ref{ex!ex2}, but specialized to contain
$\PP^2_{\Span{x_1,x_2,x_3}}$. Its equations are the maximal Pfaffians of the
Jerry$_{12}$ matrix \eqref{eq!aPf}. One checks that $Y'$ has 4 ordinary
nodes on $\PP^2$ as its only singularities for general choices of
$(a_2,a_3,b_1,b_2)$, so that it unprojects to a quasismooth $X$. Viewed
from $X$, this starts from the ordinary blowup of a general point.
\end{enumerate}

\end{document}